# Online data processing: comparison of Bayesian regularized particle filters

Roberto Casarin
Department of Economics, University of Brescia

Jean-Michel Marin*
INRIA Futurs, Projet select, Université Paris-Sud

October 25, 2018

## Abstract

The aim of this paper is to compare three regularized particle filters in an online data processing context. We carry out the comparison in terms of hidden states filtering and parameters estimation, considering a Bayesian paradigm and a univariate stochastic volatility model. We discuss the use of an improper prior distribution in the initialization of the filtering procedure and show that the regularized Auxiliary Particle Filter (APF) outperforms the regularized Sequential Importance Sampling (SIS) and the regularized Sampling Importance Resampling (SIR).

## 1 Introduction

The analysis of phenomena, which evolve over time is a common problem to many fields like engineering, physics, biology, statistics, economics and finance. A time varying system can be represented through a dynamic model, which is constituted by an observable component and an unobservable internal state. The hidden states (or latent variables) represent the informations we want to extrapolate from the observations.

In time series analysis, many approaches have been used for the estimation of dynamics models. The seminal works of Kalman (1960) and Kalman and Bucy (1960) introduce filtering techniques (the Kalman-Bucy filter) for continuous valued, linear and Gaussian dynamic systems. Maybeck (1982) motivates the use of stochastic dynamic systems in engineering and examines the estimation problems for state space models, in both a continuous and a discrete time framework. In economics, Harvey (1989) studies the state space representation of dynamic structural models and uses Kalman filter for hidden states filtering. Hamilton (1989) analyzes nonlinear time series models and introduces a filter (Hamilton-Kitagawa filter) for discrete time and discrete valued dynamic systems with a finite number of states.

In this paper, the online data processing problem is considered. In these situations, as pointed out by Liu and Chen (1998), Markov Chain Monte Carlo (MCMC) samplers are too much time demanding. To overcome this difficulty, some sequential Monte Carlo techniques have been recently developed. Doucet et al. (2001) provide the state of the art on these methods. They discuss both applications and theoretical convergence of the algorithms.

*Corresponding author: Université d'Orsay, Laboratoire de Mathématiques (Bât. 425), 91405 Orsay Cedex
jean-michel.marin@inria.fr

The contribution of this work is the comparison of three types of regularized particle filters - the regularized Sequential Importance Sampling (SIS), the regularized Sampling Importance Resampling (SIR) and the regularized Auxiliary Particle Filter (APF) - when the model parameters are unknown. The online estimation of model parameters is a difficult task (Kitagawa (1998); Storvik (2002); Berzuini and Gilks (2001); Fearnhead (2002); Djuric et al. (2002); Storvik (2002); Andrieu and Doucet (2003); Doucet and Tadic (2003); Polson et al. (2002)). We consider here the Bayesian paradigm and the regularization (see Chen and Haykin (2002)) approach of Oudjane (2000); Liu and West (2001); Musso et al. (2001); Rossi (2004) based on a kernel approximation in the parameter-augmented state space. We also discuss the initialization of the filtering procedure.

This work is structured as follow. Section 2 introduces the general representation of a Bayesian dynamic model and presents the stochastic volatility model. Section 3 reviews some regularized particle filters. Finally, Section 4 presents the results.

## 2 Bayesian Dynamic Models

We introduce the general formulation of a Bayesian dynamic model and show some fundamental relations for Bayesian inference on it. Our definition of dynamic model is general enough to include the models analyzed in Kalman (1960), Hamilton (1994), Carter and Kohn (1994), Harrison and West (1989) and in Doucet et al. (2001). Throughout this work, we use a notation similar to that one commonly used in particle filters literature (see Doucet et al. (2001)).

We denote by $\{\mathbf{x}_t; t \in \mathbb{N}\}$, $\mathbf{x}_t \in \mathcal{X} \subseteq \mathbb{R}^{n_x}$, the hidden states of the system, by $\{\mathbf{y}_t; t \in \mathbb{N}_0\}$, $\mathbf{y}_t \in \mathcal{Y} \subseteq \mathbb{R}^{n_y}$, the observable variables and by $\{\boldsymbol{\theta}_t; t \in \mathbb{N}\}$, $\boldsymbol{\theta}_t \in \Theta \subseteq \mathbb{R}^{n_\theta}$, the parameters of the model. We denote by $\mathbf{x}_{0:t} = (\mathbf{x}_0, \ldots, \mathbf{x}_t)$ the collection of hidden states up to time $t$ and with $\mathbf{x}_{-t} = (\mathbf{x}_0, \ldots, \mathbf{x}_{t-1}, \mathbf{x}_{t+1}, \ldots, \mathbf{x}_T)$ the collection of all hidden states without the $t$-th element. We use the same notations for the observable variables and parameters.

The Bayesian state space representation of a dynamic model is given by:

$$
\begin{aligned}
\mathbf{y}_t &\sim p(\mathbf{y}_t | \mathbf{x}_t, \boldsymbol{\theta}_t, \mathbf{y}_{1:t-1}) & \text{measurement density}, \\
(\mathbf{x}_t, \boldsymbol{\theta}_t) &\sim p(\mathbf{x}_t, \boldsymbol{\theta}_t | \mathbf{x}_{0:t-1}, \boldsymbol{\theta}_{0:t-1}, \mathbf{y}_{1:t-1}) & \text{transition density}, \\
\mathbf{x}_0 &\sim p(\mathbf{x}_0 | \boldsymbol{\theta}_0) & \text{initial density}, \\
\boldsymbol{\theta}_0 &\sim \pi(\boldsymbol{\theta}_0) & \text{prior density},
\end{aligned}
$$

for $t = 1, \ldots, T$.

We suppose that $p(\mathbf{x}_t, \boldsymbol{\theta}_t | \mathbf{x}_{0:t-1}, \boldsymbol{\theta}_{0:t-1}, \mathbf{y}_{1:t-1}) = p(\mathbf{x}_t, \boldsymbol{\theta}_t | \mathbf{x}_{t-1}, \boldsymbol{\theta}_{t-1}, \mathbf{y}_{1:t-1})$. We also assume that the parameters are constant over time: the transition density of the parameters is then $\delta_{\boldsymbol{\theta}_{t-1}}(\boldsymbol{\theta}_t)$ with initial value $\boldsymbol{\theta}_0 = \boldsymbol{\theta}$, $\delta_x(y)$ denotes the Dirac's mass centered in $x$. In that case, the joint transition of hidden states and parameters is:

$$p(\mathbf{x}_t, \boldsymbol{\theta}_t | \mathbf{x}_{t-1}, \boldsymbol{\theta}_{t-1}, \mathbf{y}_{1:t-1}) = p(\mathbf{x}_t | \mathbf{x}_{t-1}, \boldsymbol{\theta}_t, \mathbf{y}_{1:t-1}) \delta_{\boldsymbol{\theta}_{t-1}}(\boldsymbol{\theta}_t).$$

Let us denote by $\mathbf{z}_t = (\mathbf{x}_t, \boldsymbol{\theta}_t)$ the parameter-augmented state vector and by $\mathcal{Z}$ the corresponding augmented state space. For such models, we are interested in the prediction, filtering and smoothing densities which are given by:

$$
\begin{align}
p(\mathbf{z}_{t+1}|\mathbf{y}_{1:t}) &= \int_{\mathcal{Z}} p(\mathbf{x}_{t+1}|\mathbf{x}_t, \boldsymbol{\theta}_{t+1}, \mathbf{y}_{1:t}) \delta_{\boldsymbol{\theta}_t}(\boldsymbol{\theta}_{t+1}) p(\mathbf{z}_t|\mathbf{y}_{1:t}) d\mathbf{z}_t \,, \tag{1} \\
p(\mathbf{y}_{t+1}|\mathbf{y}_{1:t}) &= \int_{\mathcal{Z}} p(\mathbf{y}_{t+1}|\mathbf{z}_{t+1}, \mathbf{y}_{1:t}) p(\mathbf{z}_{t+1}|\mathbf{y}_{1:t}) d\mathbf{z}_{t+1} \,, \\
p(\mathbf{z}_{t+1}|\mathbf{y}_{1:t+1}) &= \frac{p(\mathbf{y}_{t+1}|\mathbf{z}_{t+1}, \mathbf{y}_{1:t}) p(\mathbf{z}_{t+1}|\mathbf{y}_{1:t})}{p(\mathbf{y}_{t+1}|\mathbf{y}_{1:t})} \,, \tag{2} \\
p(\mathbf{z}_s|\mathbf{y}_{1:t}) &= p(\mathbf{z}_s|\mathbf{y}_{1:s}) \int_{\mathcal{Z}} \frac{p(\mathbf{z}_{s+1}|\mathbf{z}_s, \mathbf{y}_{1:s}) p(\mathbf{z}_{s+1}|\mathbf{y}_{1:t})}{p(\mathbf{z}_{s+1}|\mathbf{y}_{1:s})} d\mathbf{z}_{s+1}, \quad s < t \,.
\end{align}
$$

Due to the high number of integrals that must be solved, previous densities may be difficult to evaluate with general dynamics. Some Monte Carlo simulation methods, such as particle filters, allow us to overcome these difficulties.

As an example, let us consider the stochastic volatility model. Two of the main features of the financial time series are time varying volatility and clustering phenomena in volatility. Stochastic volatility models widely used in finance have been introduced, in order to account for these features. Let $y_t$ be the observable variable with time varying volatility and $x_t$ the stochastic log-volatility process. An example of stochastic volatility model is:

$$
\begin{align}
y_t|x_t &\sim \mathcal{N}(0, e^{x_t}) \\
x_t|x_{t-1}, \boldsymbol{\theta} &\sim \mathcal{N}(\alpha + \phi x_{t-1}, \sigma^2) \\
x_0|\boldsymbol{\theta} &\sim \mathcal{N}(0, \sigma^2/(1-\phi^2)) \\
\boldsymbol{\theta} &\sim \pi(\boldsymbol{\theta})
\end{align}
$$

where $\boldsymbol{\theta} = (\alpha, \log((1+\phi)/(1-\phi)), \log(\sigma^2))$. The choice of $\pi(\boldsymbol{\theta})$ will be discussed in Section 4. Fig. 1 shows two simulated paths of $y_t$ and $x_t$.

In the next section, we deal with the problem of parameters and states joint estimation in a kernel-regularized sequential Monte Carlo framework.

## 3 Regularized particular filters

For making inference on the Bayesian dynamic model given in Section 2 in an online data processing context, MCMC algorithms are too much time demanding. Sequential importance sampling and more advanced sequential Monte Carlo algorithms called Particle Filters (Doucet et al., 2001) represent a promising alternative. The main advantage in using particle filters is that they can deal with nonlinear models and non-Gaussian innovations. In contrast to Hidden Markov Model filters, which work on a state space discretized to a fixed grid, particle filters focus sequentially on the higher density regions of the state space. This feature is common to one of the early sequential methods, the Adaptive Importance Sampling algorithm due to West (1992, 1993).

Different particle filters exist in the literature and different simulation approaches like rejection sampling, MCMC and importance sampling, can be used for the construction of a particle filter. In this work, we present some kernel-regularized particle filters, which combine the importance sampling reasoning with a suitable modification of the importance weights. The regularization approach we use is the same than the one of Liu and West (2001) and Musso et al. (2001). This approach relies upon a kernel-based reconstruction of the empirical filtering densities which produces a systematic modification of the true importance weights.

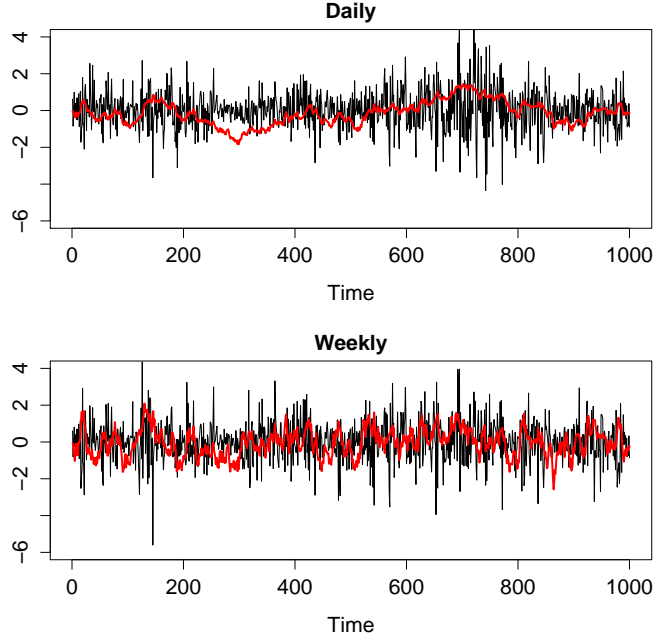

Figure 1: Simulation paths for $x_t$ (*grey line*) and $y_t$ (*black line*). Upper plot: daily dataset ($\alpha = 0$, $\phi = 0.99$ and $\sigma^2 = 0.01$). Bottom plot: weekly dataset ($\alpha = 0$, $\phi = 0.9$ and $\sigma^2 = 0.1$).

## 3.1 Regularized SIS

Let us start from the non-regularized SIS. We assume that at iteration $t > 0$ a properly weighted particle set $\{\mathbf{x}_t^i, \boldsymbol{\theta}_t^i, \gamma_t^i\}_{i=1}^N$, approximating the filtering density $p(\mathbf{x}_t, \boldsymbol{\theta}_t|\mathbf{y}_{1:t})$, is available. The empirical distribution corresponding to this approximation is:

$$p_N(\mathbf{x}_t, \boldsymbol{\theta}_t|\mathbf{y}_{1:t}) = \sum_{i=1}^N \gamma_t^i \delta_{(\mathbf{x}_t^i, \boldsymbol{\theta}_t^i)}(\mathbf{x}_t, \boldsymbol{\theta}_t). \quad (3)$$

The particles set, $\{\mathbf{x}_t^i, \boldsymbol{\theta}_t^i, \gamma_t^i\}_{i=1}^N$, can be viewed as a random discretisation of the state space $\mathcal{X} \times \Theta$ with associated probability weights $\{\gamma_t^i\}_{i=1}^N$. Thanks to this discretisation, it is possible to approximate the prediction and filtering densities given in (1) and (2):

$$p_N(\mathbf{x}_{t+1}, \boldsymbol{\theta}_{t+1}|\mathbf{y}_{1:t}) = \sum_{i=1}^N \gamma_t^i p(\mathbf{x}_{t+1}|\mathbf{x}_t^i, \boldsymbol{\theta}_{t+1}, \mathbf{y}_{1:t}) \delta_{\boldsymbol{\theta}_t^i}(\boldsymbol{\theta}_{t+1}),$$

$$p_N(\mathbf{x}_{t+1}, \boldsymbol{\theta}_{t+1}|\mathbf{y}_{1:t+1}) \propto \sum_{i=1}^N \gamma_t^i p(\mathbf{y}_{t+1}|\mathbf{x}_{t+1}, \boldsymbol{\theta}_{t+1}, \mathbf{y}_{1:t}) p(\mathbf{x}_{t+1}|\mathbf{x}_t^i, \boldsymbol{\theta}_{t+1}, \mathbf{y}_{1:t}) \delta_{\boldsymbol{\theta}_t^i}(\boldsymbol{\theta}_{t+1}).$$

The goal is now to obtain $N$ particles $\{\mathbf{x}_{t+1}^i, \boldsymbol{\theta}_{t+1}^i, \gamma_{t+1}^i\}_{i=1}^N$ from the filtering density in (2). It is proposed to sample $(\mathbf{x}_{t+1}^i, \boldsymbol{\theta}_{t+1}^i)$ according to the importance density $q(\cdot|\mathbf{x}_t^i, \boldsymbol{\theta}_t^i, \mathbf{y}_{1:t+1})$. The importance weight of particle $(\mathbf{x}_{t+1}^i, \boldsymbol{\theta}_{t+1}^i)$ is then calculated using the recursive formula:

$$\gamma_{t+1}^i \propto \gamma_t^i \frac{p(\mathbf{y}_{t+1}|\mathbf{x}_{t+1}^i, \boldsymbol{\theta}_{t+1}^i, \mathbf{y}_{1:t}) p(\mathbf{x}_{t+1}^i|\mathbf{x}_t^i, \boldsymbol{\theta}_{t+1}^i, \mathbf{y}_{1:t}) \delta_{\boldsymbol{\theta}_t^i}(\boldsymbol{\theta}_{t+1}^i)}{q(\mathbf{x}_{t+1}^i, \boldsymbol{\theta}_{t+1}^i|\mathbf{x}_t^i, \boldsymbol{\theta}_t^i, \mathbf{y}_{1:t+1})}. \quad (4)$$

The choice of an optimal importance density $q(\cdot|\mathbf{x}_t^i, \boldsymbol{\theta}_t^i, \mathbf{y}_{1:t+1})$, that is, a density which minimizes the variance of the importance weights is discussed in Pitt and Shephard (1999) and Crisan and Doucet (2000). In many cases, it is not possible to use this optimal importance density as the weight updating associated to the this density does not admit a closed-form expression. In that case, the transition density of the parameter-augmented state vector represents a natural alternative for the importance density. Indeed, the transition density represents a sort of prior at time $t$ for the parameter-augmented state vector $(\mathbf{x}_{t+1}^i, \boldsymbol{\theta}_{t+1}^i)$.

In our case, due to the presence of the Dirac point mass at the numerator of the weights it is impossible to modify over the filtering iterations the particle values for the parameters. In practice due to the loss of particle diversity in the parameter space, the weights will tend to zeros and of course stay zero for ever, so we are facing a problem of degeneracy of the empirical filtering distribution. This scenario motivates particle filtering methods known as regularized particle filters. In order to avoid the degeneracy problem and to force the exploration of the parameter space toward regions which are not covered by the prior distribution, Liu and West (2001) and Musso et al. (2001) propose to use a regularized version of the filtering density. This approach results in the modification of the weights in (4) and the definition of a new set of weights:

$$\omega_{t+1}^i \propto \omega_t^i \frac{p(\mathbf{y}_{t+1}|\mathbf{x}_{t+1}^i, \boldsymbol{\theta}_{t+1}^i, \mathbf{y}_{1:t}) p(\mathbf{x}_{t+1}^i|\mathbf{x}_t^i, \boldsymbol{\theta}_{t+1}^i, \mathbf{y}_{1:t}) K_h(\boldsymbol{\theta}_{t+1}^i - \boldsymbol{\theta}_t^i)}{q(\mathbf{x}_{t+1}^i, \boldsymbol{\theta}_{t+1}^i|\mathbf{x}_t^i, \boldsymbol{\theta}_t^i, \mathbf{y}_{1:t+1})}$$

where $K_h(y) = h^{-d} K(y/h)$ is a regularization kernel, $K$ being a positive function defined on $\mathbb{R}^{n_\theta}$ and $h$ a positive smoothing factor (bandwidth).

The modification of the importance weights defined in (4) results from two steps. The first one is the regularization of the empirical density in (3) by a kernel estimator:

$$p_N^R(\mathbf{x}_t, \boldsymbol{\theta}_t|\mathbf{y}_{1:t}) = \sum_{i=1}^N \omega_t^i \delta_{\mathbf{x}_t^i}(\mathbf{x}_t) K_h(\boldsymbol{\theta}_t - \boldsymbol{\theta}_t^i).$$

The second one is the application of an importance sampling argument to the approximated filtering density:

$$p_N^R(\mathbf{x}_{t+1}, \boldsymbol{\theta}_{t+1}|\mathbf{y}_{1:t+1}) = \sum_{i=1}^N \omega_t^i p(\mathbf{y}_{t+1}|\mathbf{x}_{t+1}, \boldsymbol{\theta}_{t+1}, \mathbf{y}_{1:t}) p(\mathbf{x}_{t+1}|\mathbf{x}_t^i, \boldsymbol{\theta}_{t+1}, \mathbf{y}_{1:t}) K_h(\boldsymbol{\theta}_{t+1} - \boldsymbol{\theta}_t^i).$$

The convergences results associated with this type of approximation are recalled in Musso et al. (2001) and Oudjane (2000). Under some usual conditions on the kernel, when the number of particles increases to infinity, the regularized empirical density converges to the right one for various criteria. For instance, we have $p_N^R \longrightarrow_{\mathcal{L}^2} p$.

Thanks to this approximation, the regularization kernel becomes the natural choice for the parameters proposal distribution. Thus, we sample $(\mathbf{x}_{t+1}^i, \boldsymbol{\theta}_{t+1}^i)$ according to:

$$q(\mathbf{x}_{t+1}|\mathbf{x}_t^i, \boldsymbol{\theta}_{t+1}, \mathbf{y}_{1:t+1}) K_h(\boldsymbol{\theta}_{t+1} - \boldsymbol{\theta}_t^i).$$

In that case, we have:

$$\omega_{t+1}^i \propto \omega_t^i \frac{p(\mathbf{y}_{t+1}|\mathbf{x}_{t+1}^i, \boldsymbol{\theta}_{t+1}^i, \mathbf{y}_{1:t}) p(\mathbf{x}_{t+1}^i|\mathbf{x}_t^i, \boldsymbol{\theta}_{t+1}^i, \mathbf{y}_{1:t})}{q(\mathbf{x}_{t+1}^i|\mathbf{x}_t^i, \boldsymbol{\theta}_{t+1}^i, \mathbf{y}_{1:t+1})}.$$

In Algorithm 1, we give a pseudo-code representation of this method.

> **Algorithm 1 - *Regularized SIS Particle Filter* -**
> · At time $t = 0$, for $i = 1, \ldots, N$, simulate $\mathbf{z}_0^i \sim p(\mathbf{z}_0)$ and set $\omega_0^i = 1/N$
> · At time $t > 0$, given $\{\mathbf{x}_t^i, \boldsymbol{\theta}_t^i, \omega_t^i\}_{i=1}^N$, for $i = 1, \ldots, N$:
>
> 1. Simulate $\boldsymbol{\theta}_{t+1}^i \sim K_h(\boldsymbol{\theta}_{t+1} - \boldsymbol{\theta}_t^i)$
>
> 2. Simulate $\mathbf{x}_{t+1}^i \sim q(\mathbf{x}_{t+1}|\mathbf{x}_t^i, \boldsymbol{\theta}_{t+1}^i, \mathbf{y}_{1:t+1})$
>
> 3. Update the weights: $\omega_{t+1}^i \propto \omega_t^i \dfrac{p(\mathbf{y}_{t+1}|\mathbf{x}_{t+1}^i, \boldsymbol{\theta}_{t+1}^i, \mathbf{y}_{1:t})p(\mathbf{x}_{t+1}^i|\mathbf{x}_t^i, \boldsymbol{\theta}_{t+1}^i, \mathbf{y}_{1:t})}{q(\mathbf{x}_{t+1}^i|\mathbf{x}_t^i, \boldsymbol{\theta}_{t+1}^i, \mathbf{y}_{1:t+1})}$.

## 3.2 Regularized SIR

As it is well known in the literature (see for example Arulampalam et al. (2001)), basic SIS algorithms have a degeneracy problem. After some iterations the empirical distribution degenerates into a Dirac's mass on a single particle. This due to the fact that the variance of the importance weights is non-decreasing over time (see Doucet et al. (2000)). In order to solve this degeneracy problem, Gordon et al. (1993) introduce the SIR algorithm. This algorithm belongs to a wider class of bootstrap filters. At each iteration, a resampling step is used to generate a new set of particles. After this resampling step, the weights of the resampled particles are uniformly distributed over the particle set.

In the initial SIR, the resampling step is done at each iteration of the algorithm. This systematic resampling can introduce extra Monte Carlo variations, see Liu and Chen (1998). This can be reduced be doing resampling only when the Effective Sample Size (ESS) is small. The ESS measures the overall efficiency of an importance sampling algorithm. The ESS is a function of the coefficient of variation of the importance weights. At iteration $t$, the empirical EES is

$$\text{ESS}_t = \frac{N}{1 + N \sum_{i=1}^N \left(\omega_t^i - N^{-1} \sum_{i=1}^N \omega_t^i\right)^2 \Big/ \left(\sum_{i=1}^N \omega_t^i\right)^2}.$$

In Algorithm 2, we give a pseudo-code representation of this method.

> **Algorithm 2 - *Regularized SIR Particle Filter* -**
> · At time $t = 0$, for $i = 1, \ldots, N$, simulate $\mathbf{z}_0^i \sim p(\mathbf{z}_0)$ and set $\omega_0^i = 1/N$
> · At time $t > 0$, given $\{\mathbf{x}_t^i, \boldsymbol{\theta}_t^i, \omega_t^i\}_{i=1}^N$, for $i = 1, \ldots, N$:
>
> 1. Simulate $\boldsymbol{\theta}_{t+1}^i \sim K_h(\boldsymbol{\theta}_{t+1} - \boldsymbol{\theta}_t^i)$
>
> 2. Simulate $\mathbf{x}_{t+1}^i \sim q(\mathbf{x}_{t+1}|\mathbf{x}_t^i, \tilde{\boldsymbol{\theta}}_{t+1}^i, \mathbf{y}_{1:t+1})$
>
> 3. Update the weights: $\omega_{t+1}^i \propto \omega_i^t \dfrac{p(\mathbf{y}_{t+1}|\tilde{\mathbf{x}}_{t+1}^i, \tilde{\boldsymbol{\theta}}_{t+1}^i, \mathbf{y}_{1:t})p(\tilde{\mathbf{x}}_{t+1}^i|\mathbf{x}_t^i, \tilde{\boldsymbol{\theta}}_{t+1}^i, \mathbf{y}_{1:t})}{q(\tilde{\mathbf{x}}_{t+1}^i|\mathbf{x}_t^i, \tilde{\boldsymbol{\theta}}_{t+1}^i, \mathbf{y}_{1:t+1})}$
>
> 4. If $ESS_{t+1} < \kappa$, simulate $\{\mathbf{x}_{t+1}^i, \boldsymbol{\theta}_{t+1}^i\}_{i=1}^N$ from $\{\mathbf{x}_{t+1}^i, \boldsymbol{\theta}_{t+1}^i, \omega_{t+1}^i\}_{i=1}^N$ (Multinomial resampling) and set $\omega_{t+1}^i = 1/N$.

The value of $\kappa < N$ is calibrated depending on the problem.

## 3.3 Regularized APF

Due to the resampling step, the basic SIR algorithm produces a progressive impoverishment (loss of diversity) of the information contained in the particle set. To overcome this difficulty, many solutions have been proposed in the literature. We refer to the APF due to Pitt and Shephard (1999) and to the regularized APF algorithm due to Liu and West (2001). In order to avoid the resampling step, the APFs use the particle index (auxiliary variable) to select most representative particles in the proposal of the new particles. The regularized joint distribution of parameter-augmented state vector and the particle index is:

$$p_N^R(\mathbf{x}_{t+1}, \boldsymbol{\theta}_{t+1}, i | \mathbf{y}_{1:t+1}) \propto p(\mathbf{y}_{t+1} | \mathbf{x}_{t+1}, \boldsymbol{\theta}_{t+1}, \mathbf{y}_{1:t}) p(\mathbf{x}_{t+1} | \mathbf{x}_t^i, \boldsymbol{\theta}_t^i, \mathbf{y}_{1:t}) K_h(\boldsymbol{\theta}_{t+1} - \boldsymbol{\theta}_t^i) \omega_t^i.$$

A sample approximating that distribution can be obtained by using the proposal:

$$q(\mathbf{x}_{t+1}^i, \boldsymbol{\theta}_{t+1}^i, j^i | \mathbf{y}_{1:t+1}) = p(\mathbf{x}_{t+1}^i | \mathbf{x}_t^{j^i}, \boldsymbol{\theta}_{t+1}^i, \mathbf{y}_{1:t}) K_h(\boldsymbol{\theta}_{t+1}^i - \boldsymbol{\theta}_t^{j^i}) q(j^i | \mathbf{y}_{1:t+1})$$

where

$$q(j^i | \mathbf{y}_{1:t+1}) \propto p(\mathbf{y}_{t+1} | \mu_{t+1}^{j^i}, m_{t+1}^{j^i}, \mathbf{y}_{1:t}) w_t^{j^i},$$

$\mu_{t+1}^{j^i}$ and $m_{t+1}^{j^i}$ are evaluated using the initial particle set. Therefore, the importance weight of particle $(\mathbf{x}_{t+1}^i, \boldsymbol{\theta}_{t+1}^i, j^i)$ is:

$$\omega_{t+1}^i \propto \frac{p(\mathbf{y}_{t+1} | \mathbf{x}_{t+1}^i, \boldsymbol{\theta}_{t+1}^i, \mathbf{y}_{1:t})}{p(\mathbf{y}_{t+1} | \mu_{t+1}^{j^i}, m_{t+1}^{j^i}, \mathbf{y}_{1:t})}.$$

In Algorithm 3 we give a pseudo-code representation of the regularized APF.

---

**Algorithm 3 - *Regularized Auxiliary Particle Filter* -**
· *At time $t = 0$, for $i = 1, \ldots, N$, simulate $\mathbf{z}_0^i \sim p(\mathbf{z}_0)$ and set $\omega_0^i = 1/N$*
· *At time $t > 0$, given $\{\mathbf{x}_t^i, \boldsymbol{\theta}_t^i, \omega_t^i\}_{i=1}^N$, for $i = 1, \ldots, N$:*

1. *Simulate $j^i \sim q(j | \mathbf{y}_{1:t+1})$ with $j \in \{1, \ldots, N\}$ (Multinomial sampling) where $\mu_{t+1}^j = \mathbb{E}(\mathbf{x}_{t+1} | \mathbf{x}_t^j, \boldsymbol{\theta}_t^j)$ and $m_{t+1}^j = \mathbb{E}(\boldsymbol{\theta}_{t+1} | \boldsymbol{\theta}_t^j)$*

2. *Simulate $\boldsymbol{\theta}_{t+1}^i \sim K_h(\boldsymbol{\theta}_{t+1} - \boldsymbol{\theta}_t^{j^i})$*

3. *Simulate $\mathbf{x}_{t+1}^i \sim p(\mathbf{x}_{t+1} | \mathbf{x}_t^{j^i}, \boldsymbol{\theta}_{t+1}^i, \mathbf{y}_{1:t})$*

4. *Update particles weights: $\omega_{t+1}^i \propto \dfrac{p(\mathbf{y}_{t+1} | \mathbf{x}_{t+1}^i, \boldsymbol{\theta}_{t+1}^i, \mathbf{y}_{1:t})}{p(\mathbf{y}_{t+1} | \mu_{t+1}^{j^i}, m_{t+1}^{j^i}, \mathbf{y}_{1:t})}.$*

---

We can say that, in the APF, the selection step is done before simulating the hidden states. This selection depends on the current value of the observable. Therefore,

- the APF is a standard way to construct a proposal distribution for the hidden states that depends on the current of the particle;

- as we will see after, to use this selection step and the transition distribution as proposal distribution for the hidden states, results in a good proposal distribution.

In the next section, we compare the performances of some regularized SIS, SIR and APF for the stochastic volatility model.

# 4 Application to the stochastic volatility model and comparison

In this section, we apply the three regularized particle filters to the stochastic volatility model presented in Section 2. We assume that the initial value of the SV process follows the stationary distribution:

$$x_0 \sim \mathcal{N}(0, \sigma^2/(1-\phi^2)).$$

For the parameters $\beta$, $\sigma$ and $\phi$, we assume the prior

$$p(\beta^2, \phi, \sigma^2) = 1/(\sigma\beta)\mathbb{I}_{(-1,1)}(\phi),$$

where $\beta = e^\alpha$. We constrain the parameter $\phi$ to take values in the open interval $(-1,1)$ in order to impose the usual stationarity condition. As we use an improper prior, it is not possible to use the prior distribution for initializing the three particle filters. We need to start with a proper weighted sample $\{x_0^i, \boldsymbol{\theta}_0^i, \omega_0^i\}_{i=1}^N$. We propose to:

1) start the sequential filtering procedure at least at the value $n$ of $t$, such that the posterior distribution of the parameters given all the observations up to time $n$ is well defined: for the considered SV model, this corresponds to set $n \geq 2$;

2) use a Markov Chain Monte Carlo (MCMC) algorithm to create a sample with uniform weights.

For the prior given above, the full conditional distributions are

$$\beta^2 | \cdots \sim \mathcal{IG}\left(\sum_{t=1}^n y_t^2 \exp(-x_t)/2, (n-1)/2\right),$$

$$\sigma^2 | \cdots \sim \mathcal{IG}\left(\sum_{t=2}^n (x_t - \phi x_{t-1})^2/2 + x_1^2(1-\phi^2), (n-1)/2\right),$$

$$\pi(\phi | \cdots) \propto (1-\phi^2)^{1/2} \exp\left(-\phi^2 \sum_{t=2}^{n-1} x_t^2 - 2\phi \sum_{t=2}^n x_t x_{t-1}\right)/2\sigma^2 \mathbb{I}_{(-1,+1)}(\phi),$$

$$\pi(x_t | \cdots) \propto \exp\left\{-\frac{1}{2\sigma^2}\left((x_t - \alpha - \phi x_{t-1})^2 - (x_{t+1} - \alpha - \phi x_t)^2\right) - \frac{1}{2}\left(x_t + y_t^2 \exp(-x_t)\right)\right\}.$$

The full conditional distributions of $\phi$ and $x_t$ are not conventional and the standard Gibbs does not apply. We propose to use the Metropolis-Hastings within Gibbs algorithm studied in Celeux et al. (2006). A detailed description of the proposal distributions for $\phi$ and $x_t$ can be found in Celeux et al. (2006). In that paper, the authors compare this MCMC scheme to an iterated importance sampling one. Note that one could alternatively use this iterated importance sampling algorithm to create a first weighted sample.

Given the initial weighted random sample $\{x_t^i, \boldsymbol{\theta}_t^i, \omega_t^i\}_{i=1}^N$, where $\boldsymbol{\theta}_t = (\alpha_t, \log((1+\phi_t)/(1-\phi_t)), \log(\sigma_t^2))$, if we use the transition density as proposal distribution for the hidden states, the regularized SIS performs the following steps:

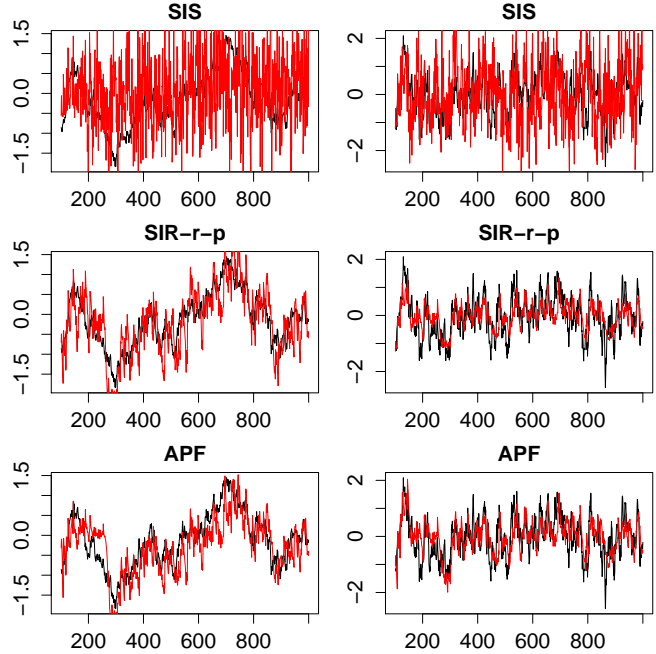

Figure 2: Daily (*left column*) and weekly (*right column*) true (*black line*) and filtered (*grey line*) log-volatility.

For $n \leq t \leq T-1$ and for $i = 1, \ldots, N$:

*(i)* Simulate $\boldsymbol{\theta}_{t+1}^i \sim \mathcal{N}\left(a\boldsymbol{\theta}_t^i + (1-a)\bar{\boldsymbol{\theta}}_t, b^2 V_t\right)$ where $V_t$ and $\bar{\boldsymbol{\theta}}_t$ are the empirical covariance matrix and the empirical mean respectively, $a \in [0,1]$ and $b^2 = (1-a^2)$,

*(ii)* Simulate $x_{t+1}^i \sim \mathcal{N}\left(\alpha_{t+1}^i + \phi_{t+1}^i x_t^i, \left(\sigma^2\right)_{t+1}^i\right)$,

*(iii)* Update the weights as follow

$$w_{t+1}^i \propto w_t^i \exp\left\{-\frac{1}{2}\left[y_{t+1}^2 \exp\left(-x_{t+1}^i\right) + x_{t+1}^i\right]\right\}.$$

In the following, we call the previous scheme SIS. Instead of the transition density, we can use a proposal distribution which depends on the current value of the observable. For instance, we can resort to the proposal distribution used for the hidden states, in the MCMC initialization step. This proposal has been introduced by Shephard and Pitt (1997). It is based on a quadratic Taylor expansion of $\exp(x_t)$, more details can be found in Shephard and Pitt (1997) or Celeux et al. (2006). It defines a new SIS which is called SIS-p in the following.

The results of a typical run of the regularized SIS on the synthetic dataset in Fig. 1, with $N = 10,000$ particles and $n = 100$ for the Gibbs initialization, are given in Fig. from 2 to 6. We can see (last row in Fig. 2) that after a few iterations the filtered log-volatility does not fit well to the true log-volatility. We measure sequentially the filtering performance of the regularized SIS by evaluating the cumulated Root Mean Square Error (RMSE). It measures the distance between the true and the filtered states and is defined as: $RMSE_t = \{\frac{1}{t}\sum_{u=1}^{t}(\hat{\mathbf{z}}_u - \mathbf{z}_u)^2\}^{\frac{1}{2}}$, where $\hat{\mathbf{z}}_t$ is the filtered state, which includes also the parameter sequential estimate. The RMSEs cumulate rapidly over time in both daily and weekly datasets (see upper and bottom plots in Fig. 3). The poor performance of the regularized SIS is due to the fact that the empirical posterior of the states and parameters degenerates into a Dirac's mass after a few iterations. The ESSs in Fig. 4 show that the regularized SIS degenerates after 30 iterations in both the daily and weekly datasets. We give some results on SIS-p in the following.

If we use the transition density as proposal distribution for the hidden states, the regularized SIR performs the following step:

For $n \leq t \leq T-1$ and for $i = 1, \ldots, N$:

(i) Simulate $\boldsymbol{\theta}_{t+1}^i \sim \mathcal{N}\left(a\boldsymbol{\theta}_t^i + (1-a)\bar{\boldsymbol{\theta}}_t, b^2 V_t\right)$ where $V_t$ and $\bar{\boldsymbol{\theta}}_t$ are the empirical covariance matrix and the empirical mean respectively and $a \in [0,1]$ and $b^2 = (1-a^2)$,

(ii) Simulate $x_{t+1}^i \sim \mathcal{N}\left(\alpha_{t+1}^i + \phi_{t+1}^i x_t^i, (\sigma^2)_{t+1}^i\right)$,

(iii) Update the weights

$$w_{t+1}^i \propto w_t^i \exp\left\{-\frac{1}{2}\left[y_{t+1}^2 \exp\{-x_{t+1}^i\} + x_{t+1}^i\right]\right\},$$

(v) If $ESS_{t+1} < \kappa$, simulate $\mathbf{z}_{t+1}^i \sim \sum_{j=1}^{N} w_{t+1}^j \delta_{\mathbf{z}_{t+1}^j}(\mathbf{z}_{t+1})$ and set $w_{t+1}^i = 1/N$.

If $\kappa = N$, the resampling step is done all the time. In that case, we call SIR the previous scheme. After some numerical experiments, we have found that a good value for $\kappa$ is $\kappa = 0.9 \times N$. In that case, the resampling step is done at regular time intervals and we called SIR-r the resulting algorithm. Moreover, as for the SIS, we can resort to the proposal of Shephard and Pitt (1997) for the hidden states. In that case, with $\kappa = N$, the corresponding algorithm called SIR-p. With $\kappa = 0.9 \times N$ the corresponding algorithm is called SIR-r-p.

The regularized APF performs the following steps:

For $n \leq t \leq T-1$ and for $i = 1, \ldots, N$:

(i) Simulate $j^i \sim q(j) \propto \sum_{k=1}^{N} w_t^k \mathcal{N}(y_{t+1}|\mu_{t+1}^k)\delta_k(j)$ where $\mu_{t+1}^k = \phi_t^k x_t^k + \alpha_t^k$,

(ii) Simulate $\boldsymbol{\theta}_{t+1}^i \sim \mathcal{N}\left(a\boldsymbol{\theta}_t^{j^i} + (1-a)\bar{\boldsymbol{\theta}}_t, b^2 V_t\right)$ where $V_t$ and $\bar{\boldsymbol{\theta}}_t$ are the empirical variance matrix and the empirical mean respectively and $a \in [0,1]$ and $b^2 = (1-a^2)$,

(iii) Simulate $x_{t+1}^i \sim \mathcal{N}\left(x_{t+1}|\alpha_{t+1}^i + \phi_{t+1}^i x_t^{j^i}, (\sigma^2)_{t+1}^i\right)$,

(iv) Update the weights

$$w_{t+1}^i \propto \exp\left\{-\frac{1}{2}\left[y_{t+1}^2\left(\exp\{-x_{t+1}^i\} - \exp\{-\mu_{t+1}^{j^i}\}\right) + x_{t+1}^i - \mu_{t+1}^{j^i}\right]\right\}.$$

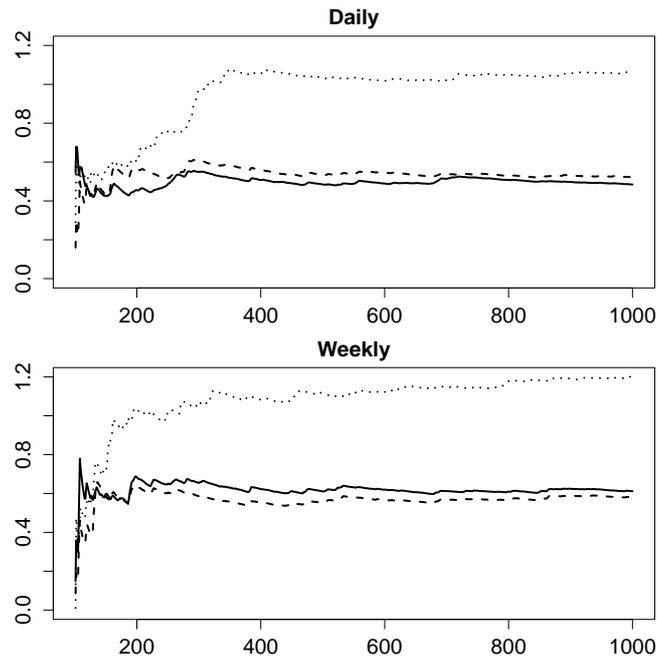

Figure 3: Daily (*upper plot*) and weekly (*bottom plot*) Root Mean Square Errors for the regularized APF (*solid line*), SIR-r-p (*dashed line*) and SIS (*dotted line*) over iterations.

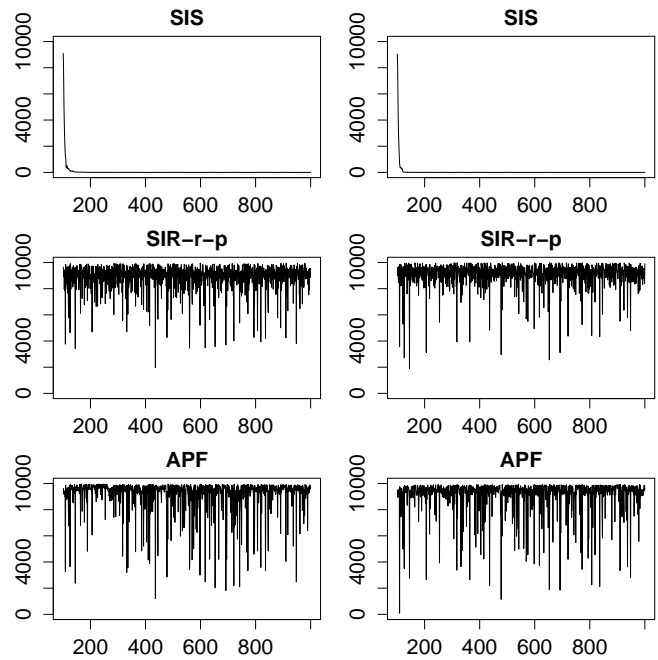

Figure 4: Daily (*left column*) and weekly (*right column*) Effective Sample Sizes over iterations.

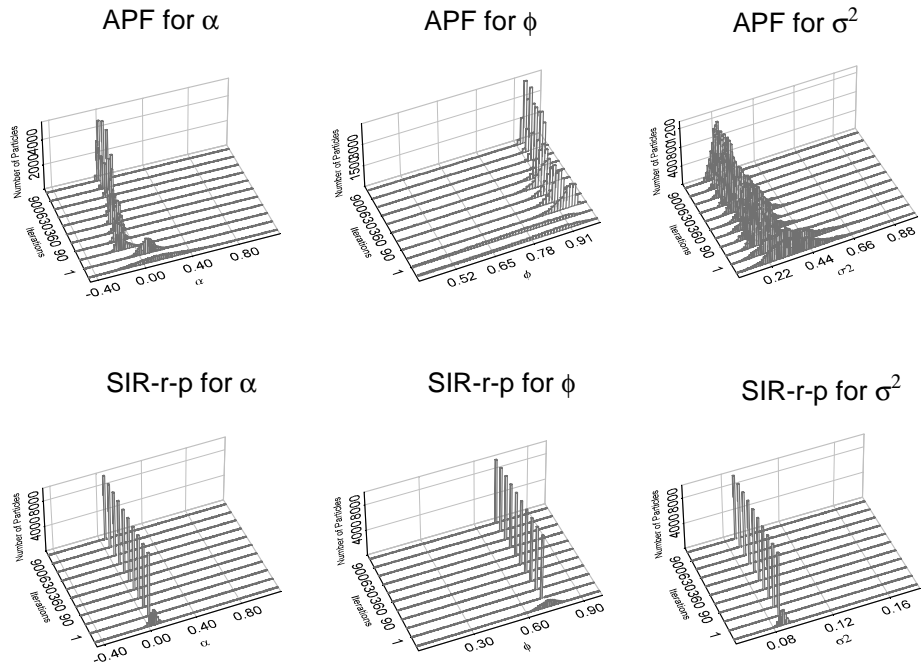

Figure 5: Evolution on daily dataset of the empirical posterior distributions of $\alpha$, $\phi$ and $\sigma^2$.

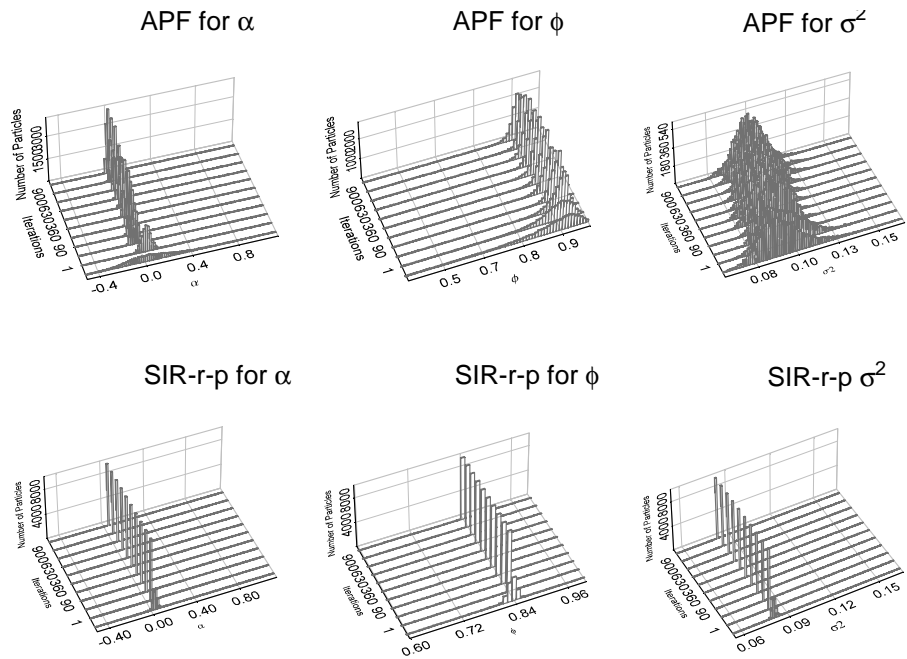

Figure 6: Evolution on weekly dataset of the empirical posterior distributions of $\alpha$, $\phi$ and $\sigma^2$.

|  | Daily Data | | | | | | |
| --- | --- | --- | --- | --- | --- | --- | --- |
| $\theta$ | SIS | SIS-p | SIR | SIR-p | SIR-r | SIR-r-p | APF |
| $\alpha$ | 0.00719 | 0.00945 | 0.00885 | 0.00925 | 0.00315 | 0.00912 | 0.00065 |
| $\phi$ | 0.66767 | 0.83264 | 0.12433 | 0.13355 | 0.15252 | 0.13456 | 0.00855 |
| $\sigma^2$ | 0.89327 | 0.87910 | 0.00676 | 0.00670 | 0.00643 | 0.00654 | 0.00506 |

Table 1: Mean Square Errors of the estimators of $\alpha$, $\phi$ and $\sigma^2$. The Mean Square Errors are estimated using the last iteration of the 10 independent runs of the filters.

|  | Weekly Data | | | | | | |
| --- | --- | --- | --- | --- | --- | --- | --- |
| $\theta$ | SIS | SIS-p | SIR | SIR-p | SIR-r | SIR-r-p | APF |
| $\alpha$ | 0.00534 | 0.00487 | 0.00589 | 0.00442 | 0.00380 | 0.00431 | 0.00016 |
| $\phi$ | 0.51290 | 0.55648 | 0.05292 | 0.03754 | 0.04885 | 0.03687 | 0.00029 |
| $\sigma^2$ | 0.70540 | 0.71242 | 0.00010 | 0.00009 | 0.00009 | 0.00009 | 0.00008 |

Table 2: Mean Square Errors of the estimators of $\alpha$, $\phi$ and $\sigma^2$. The Mean Square Errors are estimated using the last iteration of the 10 independent runs of the filters.

Note that, following Pitt and Shephard (1999), one could alternatively use in the selection step a value of $\mu_{t+1}^k$ based on the Taylor expansion of the likelihood at time $t+1$. For the three regularized particle filters, we have used a Gaussian kernel where the parameter $a$ is fixed following the usual optimal criterion.

We apply the regularized SIR-r-p and APF with $N = 10,000$ and $n = 100$ to the weekly and daily datasets of Fig. 1 and obtain the results given in Fig. from 2 to 6. The regularized SIR-r-p and APF outperform the regularized SIS in terms of ESSs and cumulated RMSEs. The ESSs can detect the degeneracy in the particle weights, but is not useful to determine the presence of another form of degeneracy, that is the absence of diversity in the particle values. The histogram of the empirical filtering distribution allows us to detect this second form of degeneracy.

As our work deals with the sequential estimation of the parameters, we choose to show the histogram of the parameters posterior. Fig. 5 and 6 exhibit the evolution over the filters iterations of the posterior of the parameters $\alpha$, $\phi$ and $\sigma^2$. In both the daily and the weekly datasets, after a few iterations the empirical posterior of the regularized SIR-r-p degenerates into a Dirac's mass.

To confirm the previous results, we have done ten independent runs of the seven algorithms: SIS, SIS-p, SIR, SIR-p, SIR-r, SIR-r-p and APF. The weekly and daily datasets vary across the 10 experiments. Our simulation study confirms the results of the single-run experiment. Fig. 7 and 8 show a comparison between the regularized schemes in terms of RMSEs. The RMSEs are estimated over the 10 independent runs of the algorithms. The regularized SIR-r-p and APF outperform the others algorithms in both the daily and the weekly datasets. The estimated Mean Square Errors for the parameters $\alpha$, $\phi$ and $\sigma^2$ (see Table 1 and 2), based on 10 independent runs of the filters, show that the regularized APF outperforms all the others schemes in term of parameters estimation.

Fig. 9 and 10 show a comparison of the filters in terms of ESS. As one could expect, in all the independent runs the regularized SIS and SIS-p weights degenerate after a few iterations. We also observe that the ESSs of the regularized SIR, SIR-r and APF are substantially equivalent. On the other hand, the average ESSs of the SIR-p and SIR-r-p are lower than the ones of the

SIR, SIR-r and APF.

## 5 Conclusion

In this work we bring into action the kernel regularization technique for particle filters and deal with the online parameter estimation problem. While the regularized APF has been already used for the parameter estimation, the regularized versions of SIS and SIR have not been considered to that aim. We focus on the joint estimation of the states and parameters and compare some algorithms on a Bayesian nonlinear model: the Bayesian stochastic volatility model. As we expected, we find evidence of the degeneracy of two different regularized SIS. Finally, we find that, in terms of parameters estimation, the regularized APF outperforms all the others schemes.

## Acknowledgments

We thank Professor Christian Robert for its helpful comments and a careful reading of a preliminary version of the paper.

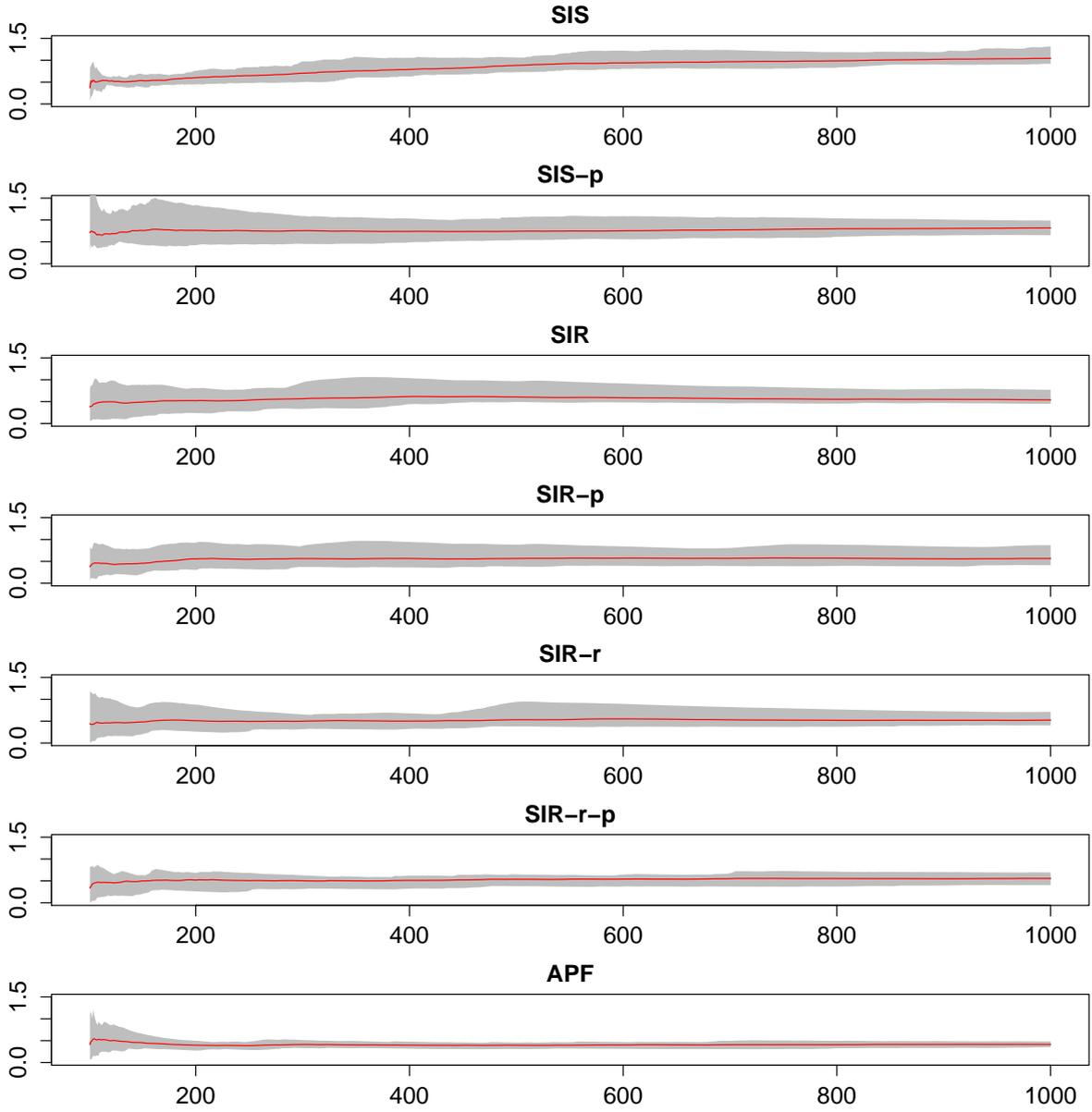

Figure 7: Comparison on daily datasets of average cumulative Root Mean Square Errors between the true and the filtered log-volatility (black line). We represent the area between maximum and minimum cumulative Root Mean Square Errors (*grey area*).

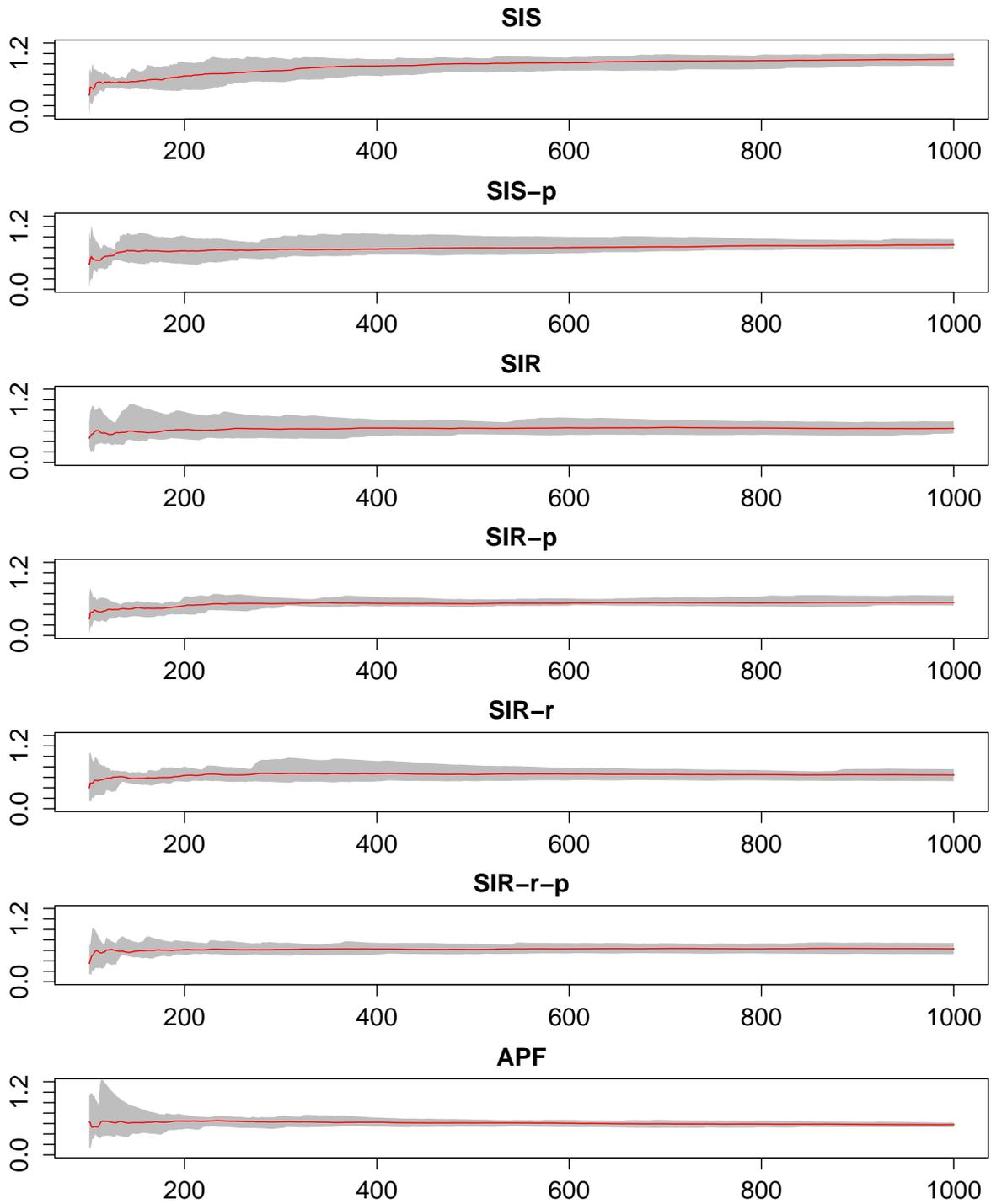

Figure 8: Comparison on weekly datasets of average cumulative RMSEs between the true and the filtered log-volatility (black line). We represent the area between maximum and minimum cumulative Root Mean Square Errors (*grey area*).

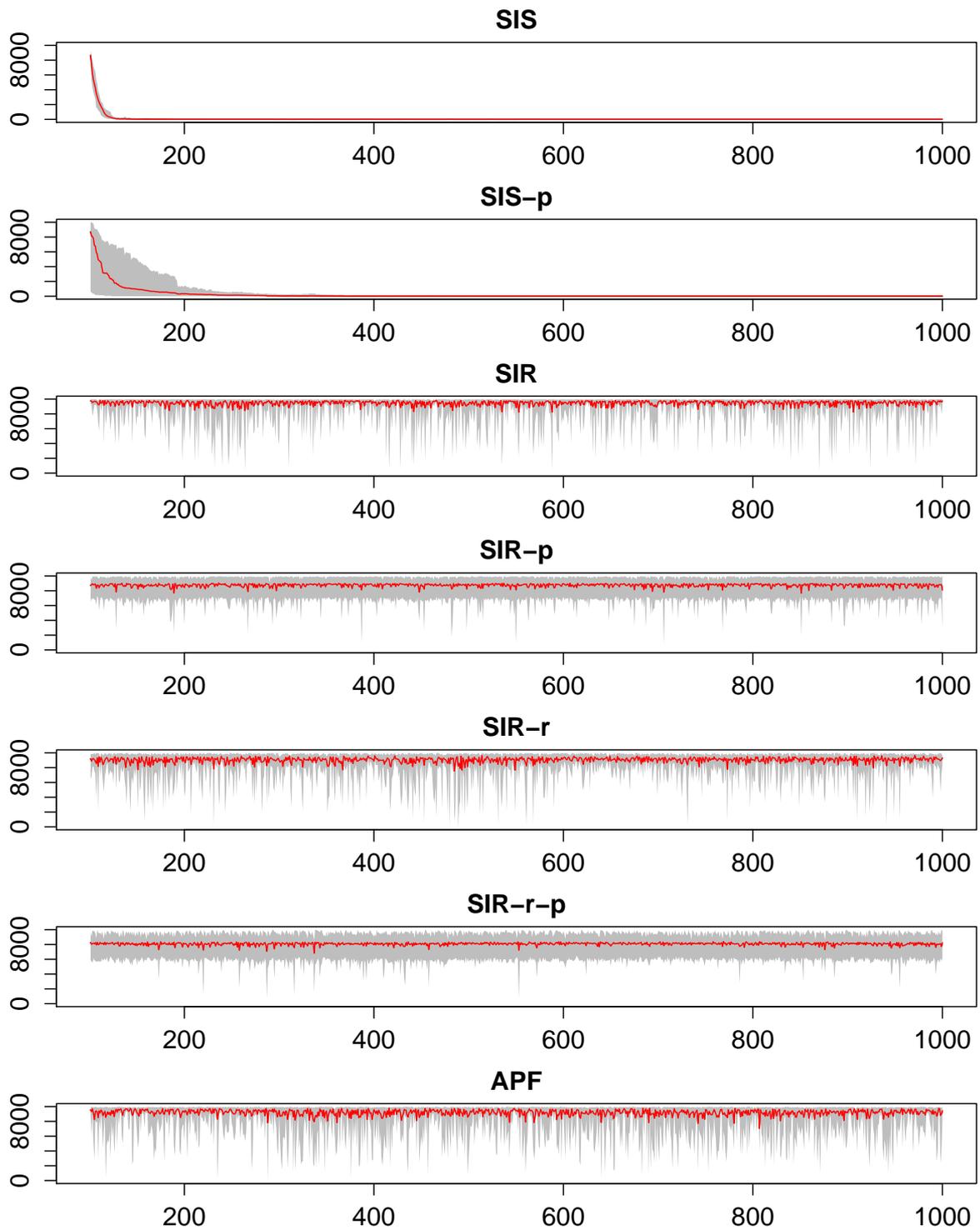

Figure 9: Comparison on daily datasets of average Effective Sample Sizes (*black line*). We represent the area between maximum and minimum Effective Sample Sizes (*grey area*).

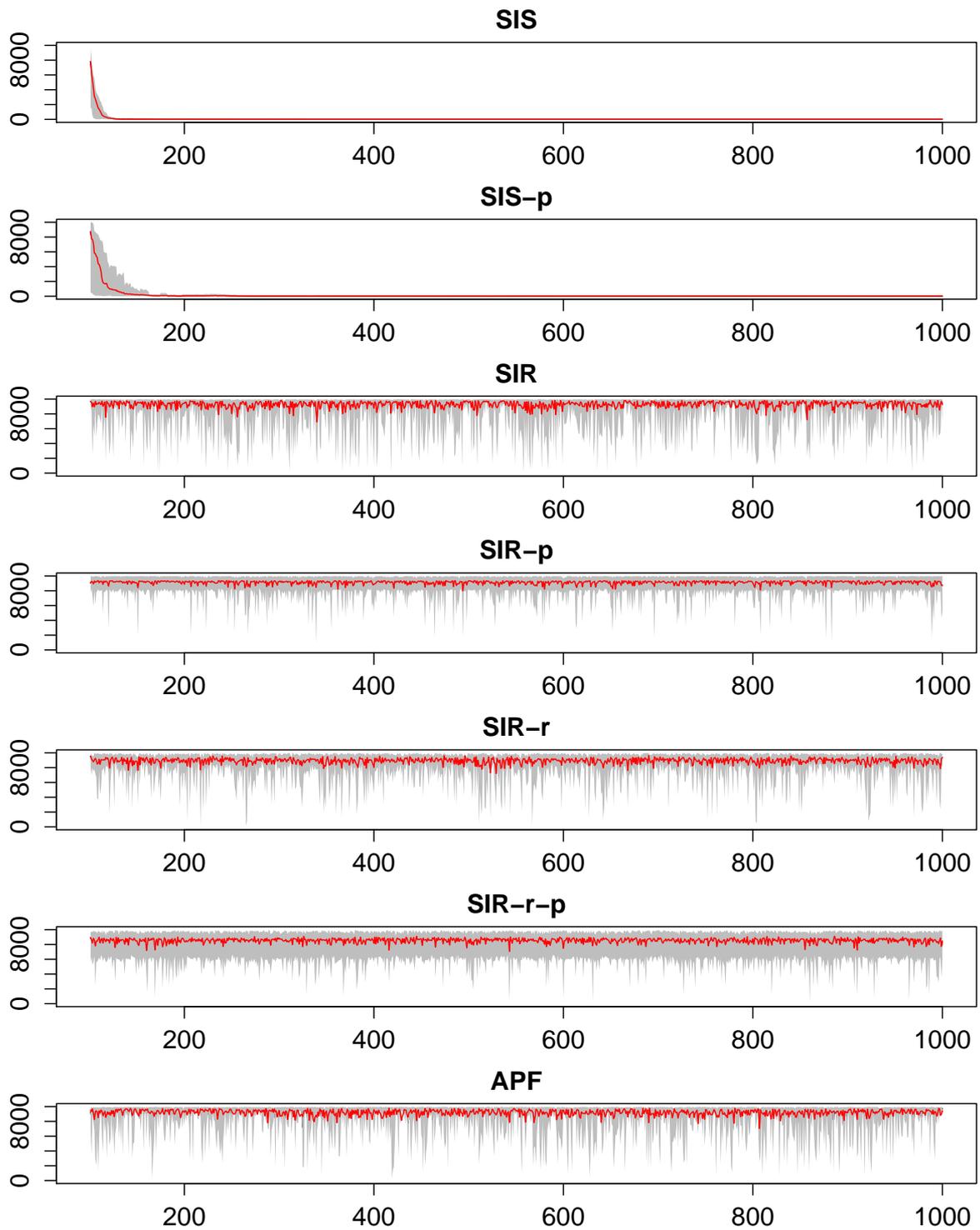

Figure 10: Comparison on weekly datasets of average Effective Sample Sizes (*black line*). We represent the area between maximum and minimum Effective Sample Size (*grey area*).

20